\documentclass{amsart}
\usepackage{amsmath}
\usepackage{amssymb}

\newtheorem{theorem}{Theorem}[section]
\newtheorem{proposition}[theorem]{Proposition}
\newtheorem{lemma}[theorem]{Lemma}
\newtheorem{corollary}[theorem]{Corollary}

\newtheorem{question}[theorem]{Question}

\theoremstyle{definition}
\newtheorem{definition}[theorem]{Definition}

\newtheorem{remark}[theorem]{Remark}

\numberwithin{equation}{section}

\newcommand{\N}{\mathbb{N}}                        
\newcommand{\C}{\mathbb{C}}                        
\newcommand{\vp}{\varphi}                          
\newcommand{\VV}{\mathcal{V}}                      
\newcommand{\FF}{\mathcal{F}}                      
\newcommand{\GG}{\mathcal{G}}                      
\newcommand{\OO}{\mathcal{O}}                      
\newcommand{\OXxi}{{\mathcal{O}}_{X,\xi}}          
\newcommand{\OYeta}{{\mathcal{O}}_{Y,\eta}}        
\newcommand{\Xpn}{X^{\{n\}}}                       
\newcommand{\Xpd}{X^{\{d\}}}                       
\newcommand{\Xxipn}{X_{{\xi}^{\{n\}}}^{\{n\}}}     
\newcommand{\Xxipd}{X_{{\xi}^{\{d\}}}^{\{d\}}}     
\newcommand{\xipn}{{\xi}^{\{n\}}}                  
\newcommand{\xipd}{{\xi}^{\{d\}}}                  
\newcommand{\vpn}{\vp^{\{n\}}}                     
\newcommand{\vpd}{\vp^{\{d\}}}                     
\newcommand{\vpnxin}{\vp^{\{n\}}_{\xi^{\{n\}}}}    
\newcommand{\vpdxid}{\vp^{\{d\}}_{\xi^{\{d\}}}}    
\newcommand{\anTor}{\widetilde{\mathrm{Tor}}}      
\newcommand{\TorR}{{\widetilde{\mathrm{Tor}}}^R}   
\newcommand{\antens}{\tilde{\otimes}}              
\newcommand{\tensR}{\tilde{\otimes}_R}             
\newcommand{\tensRn}{\tilde{\otimes}^n_R}          
\newcommand{\FpnR}{F^{\tilde{\otimes}^n_R}}        
\newcommand{\gFpnR}{(gF)^{\tilde{\otimes}^n_R}}    
\newcommand{\fd}{\mathrm{fd}}                      
\newcommand{\mm}{\mathfrak{m}}                     
\newcommand{\nn}{\mathfrak{n}}                     
\newcommand{\mb}{\mathfrak{b}}                     
\newcommand{\pp}{\mathfrak{p}}                     
\newcommand{\qq}{\mathfrak{q}}                     
\newcommand{\fbd}{\mathrm{fbd}}                    
\newcommand{\Ann}{\mathrm{Ann}}                    
\newcommand{\Spec}{\mathrm{Spec}\,}                

\begin{document}

\title{Geometric Auslander criterion for flatness}

\author{Janusz Adamus, Edward Bierstone and Pierre D. Milman}
\address{J. Adamus, Department of Mathematics, The University of Western Ontario, London, Ontario, Canada N6A 5B7
         and Institute of Mathematics, Jagiellonian University, ul. {\L}ojasiewicza 6, 30-348 Krak{\'o}w, Poland}
\email{jadamus@uwo.ca}
\address{E. Bierstone, The Fields Institute, 222 College Street, Toronto, Ontario,
Canada M5T 3J1 and Department of Mathematics, University of Toronto, Toronto, 
Ontario, Canada M5S 2E4}
\email{bierston@fields.utoronto.ca}
\address{P.D. Milman, Department of Mathematics, University of Toronto, Toronto, 
Ontario, Canada M5S 2E4}
\email{milman@math.toronto.edu}
\thanks{Research partially supported by Natural Sciences and Engineering 
Research Council of Canada Discovery Grant OGP 355418-2008 and Polish Ministry of 
Science Discovery Grant N201 020 31/1768 (Adamus), and NSERC Discovery Grants
OGP 0009070 (Bierstone), OGP 0008949 (Milman)}

\keywords{flat, torsion-free, fibred power, vertical component, analytic tensor 
product}

\begin{abstract}
Our aim is to understand the algebraic notion of flatness in explicit geometric 
terms. Let $\vp: X \to Y$ be a morphism of complex-analytic spaces, where $Y$ 
is smooth. We prove that nonflatness of $\vp$ is equivalent to a severe
discontinuity of the fibres --- the existence of a \emph{vertical component}
(a local irreducible component at a point of the source whose image is
nowhere-dense in $Y$) --- after passage to the $n$-fold fibred power of $\vp$,
where $n = \dim Y$. Our main theorem is a more general criterion for flatness
over $Y$ of a coherent sheaf of modules $\FF$ on $X$. In the case that 
$\vp$ is a morphism of complex algebraic varieties, the result implies that
the stalk $\FF_\xi$ of $\FF$ at a point $\xi \in X$ is flat over 
$R := \OO_{Y,\vp(\xi)}$ if and only if its $n$-fold tensor power is a 
torsion-free $R$-module (conjecture of Vasconcelos in the case of $\C$-algebras).
\end{abstract}
\maketitle
\setcounter{tocdepth}{1}
\tableofcontents

\section{Introduction}
\label{sec:intro}

Flatness is a subtle algebraic notion that expresses continuity
of the fibres of a mapping, and therefore the idea of a family
of varieties parametrized by a given variety $Y$.
But flatness has remained geometrically elusive because ``it
depends on infinitesimal data which are frequently invisible
at the level of topology'' (Koll\'ar \cite{Ko}). This article is
a contribution towards attempts to understand flatness in 
geometric terms. The following is a special case of our main result (see
Theorem \ref{thm:main} and Corollary \ref{cor:main}).

\begin{theorem}
\label{thm:spmain}
Let $\vp: X \to Y$ denote a morphism of complex-analytic spaces, where $Y$
is smooth, and let $\xi \in X$. Let $\vpn: \Xpn \to Y$ denote the induced
morphism from the $n$-fold fibred power of $X$ over $Y$, where $n = \dim Y$,
and let $\xipn \in \Xpn$ denote the diagonal point corresponding to $\xi$.
Then $\vp$ is not flat at $\xi$ if and only if $\vpn$ has a \emph{vertical
component} at $\xipn$; i.e., a local irreducible component (perhaps embedded)
of $\Xpn$ at $\xipn$ whose image is nowhere-dense in Y.  
\end{theorem}

Theorem \ref{thm:main} and Corollary \ref{cor:main} provide a more general
criterion for $\OO_Y$-flatness of a coherent sheaf of $\OO_X$-modules.
The theorem is an extension of a result of Galligo and Kwieci{\'n}ski \cite{GK}.

Our work is of origin in M. Auslander's criterion for freeness of a finitely generated
module over a regular local ring:

\begin{theorem}[{\cite[Thm.\,3.2]{Au}}]
\label{thm:Auslander}
Let $R$ be a regular local ring of dimension $n>0$, and let $F$ be a finite $R$-module. Then $F$ is $R$-free if and only if the $n$-fold tensor power $F^{\otimes^n_R}$ is a torsion-free $R$-module.
\end{theorem}

(Theorem~\ref{thm:Auslander} was proved in the case that $R$ is unramified by
Auslander, and extended to arbitrary regular local rings by Lichtenbaum~\cite{L}.)

\emph{Flat} is the appropriate analogue of \emph{free} for modules that 
are not necessarily finitely generated --- ``flatness $\dots$ [embodies]
that part of freeness which can be expressed in terms of linear
equations'' (Mumford \cite{Mu}). Flat is equivalent to free for
finite modules over a local ring.  
Vasconcelos \cite{V1} and Kwieci{\'n}ski~\cite{K} were the first to consider 
extending Auslander's criterion from finite modules to finite algebras 
over a regular local ring. 

We use our main theorem to prove the following result --- a generalization of 
Vasconcelos's conjecture \cite[Conj.\,6.2]{V1}, \cite[Conj.\,2.6.1]{V2} in
the case of $\C$-algebras. (See \S1.2).

\begin{theorem}
\label{thm:mainalg}
Let $R$ be a regular $\C$-algebra of finite type. 
Let $A$ denote an $R$-algebra essentially of finite type,
and let $F$ denote a finitely generated $A$-module.
Then $F$ is $R$-flat if and only if the $n$-fold tensor power
$F^{\otimes^n_R}$ is a torsion-free $R$-module, where $n = \dim R$.
\end{theorem}

An $R$-algebra \emph{essentially of finite type} means a localization 
of an $R$-algebra of finite type.

\begin{remark}
\label{rem:effective}
By Theorem~\ref{thm:mainalg} and the prime avoidance lemma \cite[Lemma\,3.3]{Eis},
in order to verify that $F$ is not $R$-flat,
it is enough to find an associated prime of
$F^{\otimes^n_R}$ in $A^{\otimes^n_R}$ which contains a 
nonzero element $r\in R$.
Thus Theorem~\ref{thm:mainalg} together with Gr\"obner-basis
algorithms for primary decomposition (see~\cite{V2} or~\cite{GP}) provides a tool
for checking flatness by effective computation.
\end{remark}

Frisch's generic flatness theorem \cite[Prop.\,VI,14]{Fri} plays an
important part in the proof of our main theorem, in the proof of 
\emph{verticality} of torsion modules \cite[Prop.\,4.5]{GK} (see 
Proposition 3.1(4) below) --- a property which is immediate in the
case of finitely generated modules over an integral domain.

The assertion of Theorem~\ref{thm:mainalg} for any field of characteristic zero follows from Theorem~\ref{thm:mainalg} as stated,
using the Tarski--Lefschetz Principle (see~\cite[Thm.\,2.1]{ABM1}). L. Avramov and S. Iyengar have more recently proved Theorem~\ref{thm:mainalg}
for an arbitrary field \cite{AI}.
In a subsequent paper~\cite{ABM2}, we prove another geometric flatness criterion using
techniques which differ from but share a common approach with those here ---
successive reduction in fibre dimension via Weierstrass preparation, to eventually
reduce the problem to a question about flatness of finitely generated modules that
is well-understood. The inductive criterion is a convenient tool to prove generic 
flatness results like Frisch's theorem.

\subsection{Vertical components in fibred powers}
Kwieci{\'n}ski~\cite{K} proved that, if $R$ is a finitely generated $\C$-algebra which is a normal domain and $A$ is a finitely generated $R$-algebra, then $A$ is $R$-flat if and only if \emph{all} tensor powers $A^{\otimes^k_R}$ are $R$-torsion-free. He used 
techniques of complex-analytic geometry, introducing the idea of a \emph{vertical
component} of an analytic mapping as a geometric analogue of torsion in commutative
algebra. Following \cite{A1}, we distinguish algebraic and geometric versions of 
Kwieci{\'n}ski's idea: Let $\vp_{\xi}:X_{\xi}\to Y_{\eta}$ denote a morphism of germs of complex-analytic spaces. 

\begin{definition}
\label{def:vert}
Let $W_{\xi}$ denote an irreducible component of $X_{\xi}$ (isolated or embedded). We say that $W_{\xi}$ is an \emph{algebraic} (respectively, \emph{geometric}) \emph{vertical component} of $\vp_{\xi}$ (or over $ Y_{\eta}$) if  $\vp_{\xi}$ maps $W_{\xi}$ to a proper analytic (respectively, nowhere-dense) subgerm of $Y_{\eta}$.
\end{definition}

(We are allowing ourselves some imprecision of language. The notation $\vp_{\xi}:X_{\xi}\to Y_{\eta}$ is meant to suggest the germ at a point $\xi \in X$ of a morphism of
complex-analytic spaces $\vp: X \to Y$, $\eta = \vp(\xi)$. In particular, we will write 
$\mathcal{O}_{X,\xi}$ for the local ring of $X_{\xi}$. Consider a representative
$\vp: X \to Y$ of $\vp_\xi$. The ``if'' clause in Definition~\ref{def:vert} means more precisely that, for a sufficiently small representative $W$ of $W_{\xi}$ in $X$, the germ
$\vp(W)_\eta$ lies in a proper analytic subgerm of $Y_\eta$ (respectively, $\vp(W)$ is nowhere-dense in $Y$.)

By the prime avoidance lemma, $\varphi_{\xi}:X_{\xi}\to Y_{\eta}$ has an algebraic
(respectively, geometric) vertical component if and only if there exists a nonzero element $m\in\mathcal{O}_{X,\xi}$ such that the zero-set germ $\VV(\Ann_{\mathcal{O}_{X,\xi}}(m))$ of the annihilator of $m$ in $\OXxi$ is mapped to a proper analytic (respectively, nowhere-dense) subgerm of $Y_{\eta}$.

We can extend the notion of vertical component to a finitely generated $\OXxi$-module
$F$:

\begin{definition}
\label{def:vert-element}
Let $I:=\Ann_{\OXxi}(F)$ and let $Z_\xi$ be the germ of a complex analytic subspace of $X$, defined by $\OO_{Z,\xi}:=\OXxi/I$.
We say that $F$ has an \emph{algebraic} (respectively, \emph{geometric}) \emph{vertical component} over $Y_{\eta}$ (or over
$\OYeta$) if $Z_{\xi}$ has an algebraic (respectively, geometric) vertical component over $Y_{\eta}$ in the sense of Definition~\ref{def:vert}; equivalently (by prime avoidance again), there exists a nonzero $m\in F$ such that the $\VV(\Ann_{\OXxi}(m))$ is mapped to a proper analytic (respectively, nowhere-dense) subgerm of $Y_{\eta}$. In the geometric case, we will call such $m$ a \emph{geometric vertical element} (or simply a \emph{vertical element}) of $F$ over $Y_{\eta}$ (or over $\OYeta$).
\end{definition}

Note that an analogous ``algebraic vertical element'' of $F$ over $Y_\eta$ is simply a (nonzero)
zero-divisor of $F$ over $\OYeta$, so there is no need to define algebraic vertical elements. A \emph{vertical element} will always mean geometric vertical.

\begin{remark}
\label{rem:vert-equivalence}
In the special case that $F = \OXxi$, $X_\xi$ has no geometric (respectively, algebraic)
vertical components over $Y_\eta$ if and only if $\OXxi$ (as an $\OXxi$-module) has no vertical elements (respectively, no zero-divisors) over $\OYeta$.
\end{remark}

Now let $R$ denote a regular local analytic $\C$-algebra of dimension $n$. Then $R$ is
isomorphic to the ring $\C\{y\} = \C\{y_1,\dots,y_n\}$ of convergent power series in $n$ variables. 
A \emph{local analytic $R$-algebra} $A$ means a ring of the form 
$R\{x\}/I=\C\{y,x\}/I$, where $I$ is an ideal in 
$\C\{y,x\} = \C\{y_1,\dots,y_n,x_1,\dots,x_m\}$, with the canonical homomorphism 
$R \to A$. Let $F$ denote an $R$-module. We say that $F$ is an \emph{almost 
finitely generated $R$-module} (following \cite{GK}) if $F$ is a 
finitely generated $A$-module, for some local analytic $R$-algebra $A$. In this 
case, there is a morphism of germs of analytic spaces 
$\vp_{\xi}:X_{\xi}\to Y_{\eta}$ such that $R \cong \OYeta$, $A \cong \OXxi$, 
$R \to A$ is the induced homomorphism $\vp^*_{\xi}: \OYeta \to \OXxi$, and $F$ 
is a finitely generated $\OXxi$-module. We say that a nonzero element $m\in F$ is \emph{vertical} over $R$ if $m$ is vertical over $\OYeta$ in the sense of Definition~\ref{def:vert-element}.

\begin{remark}
\label{rem:vert-element}
It is easy to see that the notion of vertical element is well-defined; i.e., independent of a
choice of local $R$-algebra $A$ such that $F$ is a finitely generated $A$-module.
In particular, given an almost finitely generated $R$-module $F$, we can assume 
without loss of generality that $F$ is finitely generated over the regular ring $A=R\{x\}\cong\C\{y,x\}$, where $x=(x_1,\dots,x_m)$, for some $m\geq0$.
\end{remark}

\subsection{Main results}
Our main theorem is the following flatness criterion.

\begin{theorem}
\label{thm:main}
Let $R$ be a regular local analytic $\C$-algebra and let $F$ denote an almost finitely generated $R$-module. Let $n = \dim R$. Then
$F$ is $R$-flat if and only if the $n$-fold analytic tensor power $\FpnR$ has no vertical elements over $R$.
\end{theorem}

(See Section 2 for the notion of analytic tensor power.) Theorem~\ref{thm:main} can be restated as follows.

\begin{corollary}
\label{cor:main}
Let $\vp:X\to Y$ denote a morphism of complex-analytic spaces, where $Y$ is smooth, and let $\FF$ denote a coherent sheaf of $\OO_X$-modules. Let $\xi\in X$ and $\eta=\vp(\xi)$. Then $\FF_{\xi}$ is $\OYeta$-flat if and only if the $n$-fold analytic tensor power $\FF_{\xi}^{\antens^n_{\OYeta}}$ has no vertical elements over $\OYeta$, where $n=\dim_{\eta}Y$.
\end{corollary}

Theorem~\ref{thm:spmain} is equivalent to Corollary~\ref{cor:main} in the special case that $\FF=\OO_X$,
according to Remark~\ref{rem:vert-equivalence} and the canonical isomorphism
$\OXxi\antens_{\OYeta}\dots\antens_{\OYeta}\OXxi \cong \OO_{\Xpn,\xipn}$.
\par

Theorem~\ref{thm:spmain} in the special case that $X_\xi$ is equidimensional is 
the theorem of Galligo and Kwieci{\'n}ski~\cite{GK}.  The assumption that $X_\xi$
is equidimensional guarantees that $A = \OXxi$ is a finite \emph{torsion-free} 
module over some $R$-flat algebra $S$ (where $R = \OYeta$), and Auslander's 
techniques can be extended to this case.

We are happy to acknowledge the influence of \cite{GK} on our paper.
We are also grateful to L. Avramov and S. Iyengar for pointing out
an error in an earlier version of the article (\texttt{arXiv:0901.2744v2}).
\par

We prove Theorem \ref{thm:main} in Section~\ref{sec:proof-vert-comps} below.
We reduce to the case that $A=R\{x_1,\dots,x_m\}$ using Remark~\ref{rem:vert-element}, and then argue by induction on $m$.
The case $m=0$ follows directly from Auslander's theorem.
We divide the inductive step into three cases, according as $F$ is torsion-free over $A$, $F$ is a torsion $A$-module, or neither.
The proof of the first case is independent of the inductive hypothesis, and 
again uses the argument of \cite{Au} (in the
same way it is used in \cite{GK}). The second case follows from the inductive 
assumption using Weierstrass Preparation Theorem.
\par

The most difficult new situation in the inductive step is the case that $F$ is 
neither torsion-free nor a torsion $A$-module. Our proof involves the 
geometry of the support of $F$ over $A$. 
If $F$ is $A$-torsion-free, then the support of $F$ coincides with that of 
$A$. In general, an analysis of the support of 
$F$ allows us either to reduce the fibre dimension $m$ over $R$,
or otherwise to use Proposition~\ref{prop:ness} on the variation of fibre 
dimension to produce zero-divisors over $R$ already 
for topological reasons:

We show (in Section~\ref{sec:open-vert-comps}) how nonconstancy of the dimension of the fibres 
of a morphism $\vp: X \to Y$ leads to \emph{isolated algebraic} vertical components in fibred powers of $\vp$.

By way of comparison with Corollary~\ref{cor:main}, we note that a lack of
\emph{isolated algebraic} vertical components in the $n$-fold fibred power 
(where $n$ is the target dimension) characterizes \emph{openness} of an analytic mapping
with irreducible target (see \cite{A1}, \cite{A2}). 
Proposition~\ref{prop:openness} is a simpler version of the latter result.

Frisch's generic flatness theorem asserts that if $\FF$ is a coherent sheaf of $\OO_X$-modules over a morphism $\vp: X \to Y$ of complex-analytic spaces, then $\{\xi \in X : \FF_\xi \text{ is not } \OO_{Y,\vp(\xi)}\text{--flat}\}$ has image nowhere-dense in $Y$.
Frisch's theorem is responsible for the criterion of Theorem~\ref{thm:main} in terms of \emph{geometric}
vertical elements. The existence of a vertical element in $\FpnR$ guarantees the
existence of a zero-divisor of $F^{\tilde{\otimes}^k_R}$ over $R$, for \emph{some} $k \geq n$
(by \cite{A1} and Theorem \ref{thm:main}). The following question seems to be open.

\begin{question}
With the assumptions of Theorem~\ref{thm:main}, if $F$ is not $R$-flat, does $\FpnR$
have a zero-divisor over $R$?
\end{question}

If $\vp_{\xi}:X_{\xi}\to Y_{\eta}$ and $F$ have an underlying algebraic structure as in
Theorem~\ref{thm:mainalg}, then the notions of geometric and algebraic vertical 
components coincide, so that Theorem~\ref{thm:mainalg} follows from Theorem~\ref{thm:main}:

\begin{proof}[Proof of Theorem~\ref{thm:mainalg}]

If $F$ is $R$-flat, then $F^{{\otimes}^k_R}$ is $R$-flat and therefore
$R$-torsion-free, for all $k$. 

On the other hand, suppose that $F$ is not $R$-flat.
Since $A$ is a localization of a quotient of a polynomial $R$-algebra $B = R[x_1,\dots,x_m]$
and $F$ is a finitely generated $A$-module, then $F$ is also finite over $S^{-1}B$,
for some multiplicative subset $S$ of $B$. Therefore, $F\cong S^{-1}M$,
for some finitely generated $B$-module $M$.
Flatness and torsion-freeness are both local properties;
i.e., $F$ is $R$-flat (respectively, $R$-torsion-free) if and only if $F_{\mb}$ is $R$-flat
(respectively, $R$-torsion-free), for every maximal ideal $\mb$ of $S^{-1}B$.
Since $F$ is not $R$-flat, there is a prime ideal $\pp$ in $B$ such that $\pp\cap S=\varnothing$
and $M_{\pp}$ is not $R$-flat, and it suffices to prove that
$M_{\pp}^{{\otimes}^n_R}$ is not $R$-torsion-free.

Now, the nonflatness of $M_{\pp}$ over $R$ is equivalent to that of $M_{\nn}$,
for every maximal ideal $\nn$ in $B$ containing $\pp$ (indeed, for every such $\nn$,
we have $(M_{\nn})_{\pp}\cong M_{\pp}$, and a localization of an $R$-flat $B$-module is $R$-flat).

Consider a maximal ideal $\nn$ in $B$ containing $\pp$.
We will show that $M_{\nn}^{\otimes^n_R}$ has a zero-divisor in $R$.
Let $\vp: X \to Y$ be the morphism of complex-analytic spaces associated to the
morphism $\Spec B \to \Spec R$ and let $\FF$ be the coherent sheaf of $\OO_X$-modules associated to $M$.
Let $\xi \in X$ be the point corresponding to the maximal
ideal $\nn$ of $\Spec B$. It follows from faithful flatness of completion that $\FF_\xi$
is not $\OYeta$-flat, where $\eta = \vp(\xi)$. By Theorem~\ref{thm:main}, 
$\displaystyle{\FF_{\xi}^{\antens^n_{\OYeta}}}$ has a vertical element over $\OYeta$.
Since $\vpn$ is the holomorphic map induced by the ring homomorphism $R\to B^{\otimes^n_R}$,
it follows from Chevalley's Theorem that $\displaystyle{\FF_{\xi}^{\antens^n_{\OYeta}}}$
has a zero-divisor in $\OYeta$. Hence $M_{\nn}^{\otimes^n_R}$ has a zero-divisor in $R$, as required.

Finally, let $\qq_1,\dots,\qq_s$ be the primes in $R$ whose union 
is the set of zero-divisors of $M^{\otimes^n_R}$.
It follows from the preceding paragraph that if $\nn$ is a maximal ideal containing $\pp$,
then $\nn\cap R\supset\qq_j$, for some $1\leq j\leq s$; hence $\nn\cap R\supset\qq_1\dots\qq_s$.
Since $B$ is a Jacobson ring (see \cite[Thm.\,4.19]{Eis}), $\pp$ is the intersection
of all maximal ideals containing $\pp$, and consequently $\pp\cap R\supset\qq_1\dots\qq_s$.
Then $\pp\cap R\supset\qq_j$, for some $j$, because $\pp\cap R$ is prime.
Therefore the zero-divisors from $\qq_j$ do not vanish after localizing in $\pp$,
and hence $M_{\pp}^{\otimes^n_R}$ has a zero-divisor in $R$.
\end{proof}

\section{Analytic tensor product and fibred product}
\label{sec:prod}
We briefly recall the concepts of analytic tensor product and fibred product of analytic spaces, which are used throughout the paper. 

The analytic tensor product is defined in the category of finitely generated modules over local analytic $\C$-algebras (i.e., rings of the form $\C\{z_1,\dots,z_n\}/I$ for some ideal $I$) by the usual universal mapping property for tensor product (cf.~\cite{GR}): Let $\vp_i:R\to A_i$ ($i=1,2$) be homomorphisms of local analytic $\C$-algebras. Then there is a unique (up to isomorphism) local analytic $\C$-algebra $A_1\tensR A_2$, together with homomorphisms $\theta_i:A_i\to A_1\tensR A_2$ ($i=1,2$), such that (1) $\theta_1\circ\vp_1=\theta_2\circ\vp_2$, and (2) for every pair of homomorphisms of local analytic $\C$-algebras $\psi_1:A_1\to B$, $\psi_2:A_2\to B$ satisfying $\psi_1\circ\vp_1=\psi_2\circ\vp_2$, there is a unique homomorphism of local analytic $\C$-algebras $\psi:A_1\tensR A_2\to B$ making the associated diagram commute. The algebra $A_1\tensR A_2$ is called the \emph{analytic tensor product} of $A_1$ and $A_2$ \emph{over} $R$.

For finite modules $M_1$ and $M_2$ over local analytic $R$-algebras $A_1$ and 
$A_2$, respectively, there is a unique (up to isomorphism) finite $A_1\tensR A_2$-module $M_1\tensR M_2$, together with an $R$-bilinear mapping 
$\rho:M_1\times M_2\to M_1\tensR M_2$, such that for every $R$-bilinear $\kappa:M_1\times M_2\to N$, where $N$ is a finite $A_1\tensR A_2$-module, there is a unique homomorphism of $A_1\tensR A_2$-modules $\lambda:M_1\tensR M_2\to N$ satisfying $\kappa=\lambda\circ\rho$.
The module $M_1\tensR M_2$ is called the \emph{analytic tensor product} of $M_1$ and $M_2$ \emph{over} $R$.
\par

It is sometimes convenient to express the analytic tensor product of modules over a local analytic $\C$-algebra in terms of ordinary tensor product of certain naturally associated modules: Given homomorphisms of local analytic $\C$-algebras $\vp:R\to A_i$, and finitely generated $A_i$-modules $M_i$ ($i=1,2$), the modules $M_1\tensR A_2$ and $A_1\tensR M_2$ are finitely generated over $A_1\tensR A_2$, and there is a canonical isomorphism
\[
M_1\tensR M_2\ \cong\ (M_1\tensR A_2)\otimes_{A_1\tensR A_2}(A_1\tensR M_2)\,.
\]
In particular, if $A_1=R\{x\}/I_1$ and $A_2=R\{t\}/I_2$, where $x=(x_1,\dots,x_l)$,  $t=(t_1,\dots,t_m)$ are systems of variables and $I_1 \subset R\{x\}$, $I_2 \subset R\{t\}$ are ideals, then
\begin{align*}
A_1\tensR A_2\ &\cong\ (A_1\tensR R\{t\})\otimes_{R\{x\}\tensR R\{t\}}(R\{x\}\tensR A_2)\\
&\cong\ (R\{x,t\}/I_1R\{x,t\})\otimes_{R\{x,t\}}(R\{x,t\}/I_2R\{x,t\})\\
&\cong\ R\{x,t\}/(I_1R\{x,t\}+I_2R\{x,t\})\,.
\end{align*}

The fibred product of analytic spaces is defined by a dual universal mapping property (see~\cite{F}): Let $\vp:X_i\to Y$ ($i=1,2$) denote holomorphic
mappings of complex analytic spaces. Then there exists a unique (up to isomorphism) complex analytic space $X_1\times_YX_2$, together with holomorphic maps $\pi_i:X_1\times_YX_2\to X_i$ ($i=1,2$), such that (1) $\vp_1\circ\pi_1=\vp_2\circ\pi_2$, and (2) for every pair of holomorphic maps $\psi_1:X\to X_1$, $\psi_2:X\to X_2$ satisfying $\vp_1\circ\psi_1=\vp_2\circ\psi_2$, there is a unique holomorphic map $\psi:X\to X_1\times_YX_2$ making the associated diagram commute. The space $X_1\times_YX_2$ is called the \emph{fibred product} of $X_1$ and $X_2$ \emph{over} $Y$ (more precisely, \emph{over} $\vp_1$ and $\vp_2$). There is a canonical holomorphic mapping
$\vp_1\times_Y\vp_2:X_1\times_YX_2\to Y$,
given by $\vp_1\times_Y\vp_2=\vp_i\circ\pi_i$ (where $i = 1$ or $2$).
\par

Given a holomorphic map $\vp:X\to Y$ of complex analytic spaces, with $\vp(\xi)=\eta$, let $\vp_{\xi}:X_{\xi}\to Y_{\eta}$ denote the germ of $\vp$ at $\xi$. We denote by $\vpd:\Xpd\to Y$ the canonical map from the $d$-fold fibred power of $X$ over $Y$ to $Y$, and by $\vpdxid:\Xxipd\to Y_{\eta}$ its germ at the point $\xipd:=(\xi,\dots,\xi)\in\Xpd$.\par

Suppose that $\vp_1:X_1\to Y$ and $\vp_2:X_2\to Y$ are holomorphic mappings of analytic spaces, with $\vp_1(\xi_1)=\vp_2(\xi_2)=\eta$. Then the local rings $\OO_{X_i,\xi_i}$ ($i=1,2$) are $\OYeta$-modules and, by the uniqueness of fibred product and of analytic tensor product, the local ring $\OO_{Z,(\xi_1,\xi_2)}$ of the fibred product $Z=X_1\times_YX_2$ at $(\xi_1,\xi_2)$ is canonically isomorphic to
$\OO_{X_1,\xi_1}\antens_{\OYeta}\OO_{X_2,\xi_2}$.
Therefore, given a holomorphic germ $\vp_{\xi}:X_{\xi}\to Y_{\eta}$, we will identify the $d$-fold analytic tensor power $\OXxi^{{\antens}^d_{\OYeta}}=\OXxi\antens_{\OYeta}\dots\antens_{\OYeta}\OXxi$ with the local ring of the $d$-fold fibred power $\OO_{X^{\{d\}},\xi^{\{d\}}}$, for $d\geq1$.

\section{Homological properties of almost finitely generated modules}
\label{sec:almost-fin-modules}

We first recall some homological properties of almost finitely generated modules, 
established by Galligo and Kwieci{\'n}ski \cite{GK}, that generalize the 
corresponding properties of finite modules used by Auslander \cite{Au}. 
We then generalize a lemma of Auslander \cite[Lemma 3.1]{Au} to almost
finitely generated modules (Lemma 3.3 below).

Our proof of the main theorem~\ref{thm:main} in the case that $F$ is torsion-free
over $A$ (Case (1) in Section~\ref{sec:proof-vert-comps}) follows the argument of \cite{Au}.
Auslander's main tools are two addition formulas: the Auslander-Buchsbaum 
formula (see \cite[Ch.\,VII, Prop.\,1.12]{Ku}) and additivity of projective 
dimension \cite[Cor.\,1.3]{Au}. 
We replace these by tools adapted to almost finitely generated modules:
an Auslander-Buchsbaum type formula for flat dimension (Proposition~\ref{rem:GK-facts}(2))
and additivity of flat dimension (Lemma~\ref{lem:fd-sum}).
\par

Let $R=\C\{y_1,\dots,y_n\}$ denote a regular local analytic $\C$-algebra of 
dimension $n$. Let $\tensR$ denote the analytic tensor product over $R$, 
and let $\TorR$ be the corresponding derived functor.\par

Let $F$ denote an almost finitely generated $R$-module. We define the 
\emph{flat dimension} $\fd_R(F)$ of $F$ over $R$ as the minimal length of a 
flat resolution of $F$ (i.e., a resolution by $R$-flat modules). 
It is easy to see that
\begin{equation}
\fd_R(F)=\max\{i\in\N:\, \TorR_i(F,N)\neq0\ \ \mathrm{for\ some}\ N\}\,.
\end{equation}
Indeed, if $M$ is an almost finitely generated $R$-module, then
\begin{equation}
M\ \mathrm{is\ }R\mathrm{-flat}\, \Leftrightarrow\ \TorR_1(M,R/\mm_R)=0\,,
\end{equation}
where $\mm_R$ is the maximal ideal of $R$ (cf. \cite[Prop.\,6.2]{H}). 
Let $(A,\mm_A)$ be a regular local $R$-algebra such that $F$ is a finite 
$A$-module. Then $(3.1)$ follows from $(3.2)$ applied to the kernels of a 
minimal $A$-free (hence $R$-flat) resolution
\[
\FF_*:\quad \ldots\stackrel{\alpha_{i+1}}{\longrightarrow} 
F_{i+1}\stackrel{\alpha_i}{\longrightarrow} 
F_i\stackrel{\alpha_{i-1}}{\longrightarrow}\ldots\stackrel{\alpha_1}{\longrightarrow} 
F_1\stackrel{\alpha_0}{\longrightarrow} F_0\to F
\]
of $F$. ($\FF_*$ \emph{minimal} means that $\alpha_i(F_{i+1})\subset\mm_A F_i$,
for all $i\in\N$).\par

The \emph{depth} $\mathrm{depth}_R(F)$ of $F$ as an $R$-module is defined as 
the length of a maximal $F$-sequence in $R$ (i.e., a sequence 
$a_1,\dots,a_s\in\mm_R$ such that $a_j$ is not a zero-divisor in 
$F/(a_1,\dots,a_{j-1})F$, for $j=1,\dots,s$). Since all the maximal $F$-sequences 
in $R$ have the same length, depth is well defined: As observed in 
\cite[Lemma\,2.4]{GK}, the classical proof of Northcott-Rees for finitely 
generated modules (see, e.g.,~\cite[\S VI, Prop.\,3.1]{Ku}), carries over to the 
case of almost finitely generated modules.
\medskip

\begin{proposition}
\label{rem:GK-facts}
Let $M$ and $N$ be almost finitely generated $R$-modules. 
Then the following properties hold.
\begin{enumerate}
\item \emph{Rigidity of $\TorR$} \cite[Prop.\,2.2(4)]{GK}. 
If $\anTor^R_{i_0}(M,N)=0$ for some $i_0\in\N$, then $\anTor^R_i(M,N)=0$ 
for all $i\geq i_0$.

\item \emph{Auslander--Buchsbaum-type formula} \cite[Thm.\,2.7]{GK}.
\[
\mathrm{fd}_R(M)+\mathrm{depth}_R(M)=n\,.
\]

\item \emph{Additivity of flat dimension} \cite[Prop.\,2.10]{GK}.
If $\anTor^R_i(M,N)=0$ for all $i\geq1$, then
\[
\mathrm{fd}_R(M)+\mathrm{fd}_R(N)=\mathrm{fd}_R(M\tensR N)\,.
\]

\item \emph{Verticality of $\TorR$} (cf. \cite[Prop.\,4.5]{GK}). 
For all $i\geq1$, $\anTor^R_i(M,N)$ is an almost finitely generated $R$-module, 
and every element of $\anTor^R_i(M,N)$ is vertical over $R$ (recall
Definition~\ref{def:vert-element}).
\end{enumerate}
\end{proposition}

\begin{remark}
\label{rem:how-to-strengthen}
The analytic $\anTor^R_i$ need not be torsion $R$-modules, except in the
case that $M$ and $N$ are finitely generated over $R$. (In this case,
$\anTor = \mathrm{Tor}$.) It seems to be unknown whether the $\anTor^R_i$ 
necessarily contain $R$-zero-divisors (cf. Question 1.12). 
\end{remark}

\begin{lemma}
\label{lem:fd-sum}
Let $A=R\{x\}$ denote a regular local analytic $R$-algebra, $x=(x_1,\dots,x_m)$.
Let $F$ be a finitely generated $A$-torsion-free module, and let $N$ be a module 
which is finitely generated over $B=A^{\antens^j_R}$, for some $j\geq1$.
Suppose that $F\tensR N$ has no vertical elements over $R$. Then:
\begin{enumerate}
\item $N$ has no vertical elements over $R$;
\item $\TorR_i(F,N)=0$, for all $i\geq1$;
\item $\fd_R(F)+\fd_R(N)=\fd_R(F\tensR N)$.
\end{enumerate}
\end{lemma}

\begin{proof}
To prove $(1)$, consider $N'=\{n\in N: n$ is vertical over $R\}$. It is easy to 
see that $N'$ is a $B$-submodule of $N$. Indeed, if $n,n_1,n_2\in N'$ and $b\in B$,
then
\[
\Ann_B(n_1+n_2)\supset\Ann_B(n_1)\cdot\Ann_B(n_2) \quad \mathrm{and} 
\quad \Ann_B(bn)\supset\Ann_B(n)\,;
\]
hence the zero set germ $\VV(\Ann_B(n_1+n_2))$ is mapped into the union of the
(nowhere-dense) images of $\VV(\Ann_B(n_1))$ and $\VV(\Ann_B(n_2))$, 
and $\VV(\Ann_B(bn))$ is mapped into the image of $\VV(\Ann_B(n))$. 
Therefore, we get an exact sequence of $B$-modules,
\[
0\to N'\to N\to N''=N/N'\to0\,.
\]
Tensoring with $F$ induces a long exact sequence of 
$A\tensR B\cong A^{\antens^{j+1}_R}$-modules,
\begin{align}
\ldots &\to\TorR_{i+1}(F,N')\to\TorR_{i+1}(F,N)\to\TorR_{i+1}(F,N'')\\
&\to\TorR_i(F,N')\to\ldots\to\TorR_1(F,N'')\notag\\
&\to F\tensR N'\to F\tensR N\to F\tensR N''\to0\,.\notag
\end{align}
Since every element of $N'$ is vertical over $R$, the same is true for $F\tensR N'$
(indeed, 
$\Ann_{A\tensR B}(f\tensR n)\supset 1\tensR \Ann_B(n)$ for all $f\in F,n\in N'$).
But $F\tensR N$ has no vertical elements, by assumption, so that 
$F\tensR N'\to F\tensR N$ is the zero map; hence $F\tensR N\cong F\tensR N''$. 
In particular, $F\tensR N''$ has no vertical elements over $R$.\par

Since $F$ is $A$-torsion-free, there is an injection of $F$ into a finite free 
$A$-module $L$ (obtained by composing the natural map $F\to(F^*)^*$, which is 
injective in this case, with the dual of a presentation of $F^*$). The exact sequence of 
$A$-modules $0\to F\to L\to L/F\to0$ induces a long exact sequence of 
$A\tensR B$-modules,
\begin{align*}
\ldots &\to\TorR_{i+1}(L,N'')\to\TorR_{i+1}(L/F,N'')\to\TorR_i(F,N'')\\
&\to\TorR_i(L,N'')\to\ldots\to\TorR_1(L,N'')\to\TorR_1(L/F,N'')\\
&\to F\tensR N''\to L\tensR N''\to L/F\tensR N''\to0\,.
\end{align*}

Since $L$ is a free $A$-module and therefore $R$-flat, $\TorR_i(L,N'')=0$ for 
all $i\geq1$, and we obtain isomorphisms
\begin{equation}
\TorR_{i+1}(L/F,N'')\,\cong\,\TorR_i(F,N''), \quad i\geq1\,,
\end{equation}
as well as injectivity of $\TorR_1(L/F,N'')\to F\tensR N''$. But $F\tensR N''$ has 
no vertical elements, while every element of $\TorR_1(L/F,N'')$ is vertical over 
$R$ (by Prop.~\ref{rem:GK-facts}(4)); hence $\TorR_1(L/F,N'')\to F\tensR N''$ is the zero map.
Therefore, $\ldots\stackrel{0}{\rightarrow}\TorR_1(L/F,N'')\stackrel{0}{\rightarrow}\ldots$ 
is exact; hence $\TorR_1(L/F,N'')=0$. By rigidity of $\TorR$ (Prop.~\ref{rem:GK-facts}(1)), 
$\TorR_{i+1}(L/F,N'')=0$ for all $i\geq1$, so by (3.4),
\begin{equation}
\TorR_i(F,N'')=0, \quad i\geq1\,. 
\end{equation}
In particular, $\TorR_1(F,N'')=0$, hence,
$\ldots\stackrel{0}{\rightarrow}F\tensR N'\stackrel{0}{\rightarrow}\ldots$ is 
exact, by (3.3), so that $F\tensR N'=0$. However, 
$F\tensR N'\cong(F\tensR B)\otimes_{A\tensR B}(A\tensR N')$ is an (ordinary) 
tensor product of finitely generated modules over a regular local ring $A\tensR B$,
so it is zero only if one of the factors is zero. We conclude that $A\tensR N'=0$, 
and therefore $N'=0$, by $R$-flatness of $A$. This proves assertion $(1)$.\par

Now, $\TorR_i(F,N')=0$ for all $i\geq0$; hence 
$\TorR_i(F,N)\cong\TorR_i(F,N'')$ for all $i\geq1$, by (3.3). Therefore, 
$\TorR_i(F,N)=0$ for all $i\geq1$, by (3.5), proving $(2)$.

Assertion $(3)$ follows from Proposition~\ref{rem:GK-facts}(3) and $(2)$.
\end{proof}

\section{Vertical components and variation of fibre dimension}
\label{sec:open-vert-comps}

In this section, we describe a relationship between the filtration of the
target of an analytic mapping $\vp: X \to Y$ by fibre dimension and the
isolated irreducible components of the $n$-fold fibred power $X^{\{n\}}$,
where $n = \dim Y$.

Let $\vp_{\xi}:X_{\xi}\to Y_{\eta}$ be a morphism of germs of analytic spaces, 
where $Y_{\eta}$ is irreducible and of dimension $n$. Let $Y$ be an irreducible 
representative of $Y_{\eta}$, and let $X$ be a representative of $X_{\xi}$, such 
that the components of $X$ are precisely the representatives in $X$ of the 
components of $X_{\xi}$, and $\vp(X)\subset Y$, where $\vp$ represents the 
germ $\vp_{\xi}$. Let $\fbd_x\vp$ denote the \emph{fibre dimension} $\dim_x\vp^{-1}(\vp(x))$
of $\vp$ at a point $x\in X$.\par

We will use the following notation in this section: $l:=\min\{\fbd_x\vp: x\in X\}$,
$k:=\max\{\fbd_x\vp: x\in X\}$, and $A_j:=\{x\in X: \fbd_x\vp\geq j\}$, 
$l\leq j\leq k$. Then $X=A_l\supset A_{l+1}\supset\dots\supset A_k$ and,
by upper-semicontinuity of fibre dimension (see Cartan--Remmert Theorem 
\cite[{\S}V.3.3,\,Thm.\,5]{Loj}), the $A_j$ are analytic in $X$. Define 
$B_j:=f(A_j)=\{y\in Y: \dim\vp^{-1}(y)\geq j\}$, $l\leq j\leq k$. 
Upper-semicontinuity of $\fbd_x\vp$ (as a function of $x$) implies that the 
germs $(A_j)_{\xi}$ and $(B_j)_{\eta}$ are independent of the 
choices of representatives made above.\par

Note that, except for $B_k$ (cf. proof of Proposition~\ref{prop:openness} 
below), the $B_j$ may not even be semianalytic in general. This fact is 
responsible for a complicated relationship between the algebraic vertical and 
geometric vertical components in the fibred powers of $X$ over $Y$ 
(see \cite{A2} for a detailed discussion), but will not affect our considerations
here, which rely only on the properties of $B_k$.

\begin{proposition}[{\cite[Prop.\,2.1]{A1}}]
\label{prop:filtration}
Under the assumptions above, let $\bigcup_{i\in I}W_i$ 
denote the decomposition of $(\Xpn)_{\mathrm{red}}$ into finitely many 
isolated irreducible components through $\xi^{\{ n\}}$. Then:
\begin{enumerate}
\item For each $j=l,\dots,k$, there is an index subset $I_j\subset I$ such that
\[
B_j=\bigcup_{i\in I_j}\vpn(W_i)\,.
\]
\item Let $y\in B_j$ and let $s = \dim\vp^{-1}(y)$ ($s\geq j$). If $Z$ is an 
isolated irreducible component of the fibre $(\vpn)^{-1}(y)$, of dimension $ns$, 
and $W$ is an irreducible component of $\Xpn$ containing $Z$, 
then $\vpn(W)\subset B_j$.
\end{enumerate}
\end{proposition}

\begin{proof}
For (2), fix $j\geq l+1$. (The statement is trivial for $j=l$, since
$B_l=\vp(X)$.)
Suppose that there exists $x=(x_1,\dots,x_n)\in W$ such that 
$\vp(x_1)\in Y\setminus B_j$ (and hence $\vp(x_i)\in Y\setminus B_j$, $i\leq n$). 
Then $\fbd_{x_i}\vp\leq j-1$, $i=1,\dots,n$; hence $\fbd_x\vpn\leq n(j-1)=nj-n$. 
In particular, the generic fibre dimension of $\vpn|_W$ is at most $nj-n$. 
Since $\mathrm{rank}(\vpn|_W)\leq\dim Y=n$, then $\dim W\leq(nj-n)+n=nj$ 
(see, e.g.,~\cite[V.3]{Loj}).\par

Now we have $W\supset Z$, $\dim W\leq nj$, $\dim Z=ns\geq nj$, and both $W$ and 
$Z$ are irreducible analytic sets in~$\Xpn$. This is possible only if $W=Z$;
hence $\vpn(W)=\vpn(Z)=\{ y\}\subset B_j$; a contradiction. 
Therefore $\vpn(W)\subset B_j$, completing the proof of (2).\par

Part (1) follows immediately, since if $y\in B_j$ and $Z$ is an irreducible
component of $(\vpn)^{-1}(y)$ of the highest dimension, then there exists 
an isolated irreducible component $W$ of $\Xpn$ that contains $Z$.
\end{proof}
\medskip

The following is a simplified variant of an openness criterion of 
\cite[Thm.\,2.2]{A1}, proved here under somewhat weaker assumptions.

\begin{proposition}
\label{prop:openness}
Let $\vp_{\xi}:X_{\xi}\to Y_{\eta}$ be a morphism of germs of analytic spaces. 
Suppose that $Y_{\eta}$ is irreducible, $\dim Y_{\eta} = n$, $\dim X_{\xi}=m$, 
and the maximal fibre dimension of $\vp_{\xi}$ is not generic on some
$m$-dimensional irreducible component of $X_{\xi}$. Then the $n$-fold fibred
power $\vpnxin:\Xxipn\to Y_{\eta}$ contains an isolated 
algebraic vertical component.
\end{proposition}

\begin{proof}
As above, let $\vp:X\to Y$ be a representative of $\vp_{\xi}$, 
where $Y$ is irreducible and of dimension $n$. Let $k:=\max\{\fbd_x\vp: x\in X\}$, 
$A_k:=\{x\in X: \fbd_x\vp=k\}$, and $B_k:=\vp(A_k)=\{y\in Y:\dim\vp^{-1}(y)=k\}$. 
Then the fibre dimension of $\vp$ is constant on the analytic set $A_k$. By the
Remmert Rank Theorem (see \cite[{\S}V.6,\,Thm.\,1]{Loj}), $B_k$ is locally analytic in $Y$, 
of dimension $\dim A_k-k\leq\dim X-k$. Since $\eta\in B_k$, after shrinking $Y$ if necessary, we can assume that $B_k$ 
is an analytic subset of $Y$. Therefore, by Proposition~\ref{prop:filtration}, 
it is enough to show that the analytic germ $(B_k)_{\eta}$ is a proper subgerm
of $Y_{\eta}$. Let $U$ be an isolated irreducible component of $X$, of dimension 
$m=\dim X$, and such that $k$ is not the generic fibre dimension of $\vp|_U$.
It follows that
\[
\dim Y\ \geq\ \dim U-\ \mathrm{generic}\ \fbd\vp|_U\ \geq\ m-k+1\,.
\]
Then $\dim B_k\leq m-k<\dim Y$; hence $\dim (B_k)_{\eta}<\dim Y=\dim Y_{\eta}$, 
so that $(B_k)_{\eta}\varsubsetneq Y_{\eta}$.
\end{proof}

\section{Proof of the main theorem}
\label{sec:proof-vert-comps}

Let $F$ be an almost finitely generated module over $R:=\C\{y_1,\dots,y_n\}$.
By Remark~\ref{rem:vert-element}, there exists $m\geq0$ such that $F$ is finitely generated as a module over $A=R\{x\}$, 
where $x=(x_1,\dots,x_m)$. Let $X$ and $Y$ be connected open neighbourhoods of the origins in $\C^{m+n}$ and $\C^n$ (respectively),
and let $\vp:X\to Y$ be the canonical coordinate projection.
Let $\FF$ be a coherent sheaf of $\OO_X$-modules whose stalk at the origin in $X$ equals $F$.
We can identify $R$ with $\OO_{Y,0}$ and $A$ with $\OO_{X,0}$.
Then $F$ is $R$-flat if and only if $\FF_0$ is $\OO_{Y,0}$-flat.
\par

The ``only if'' direction of Theorem~\ref{thm:main} is easy to establish (see 
\S\,\ref{subsec:proof}). Our proof of the more difficult ``if'' direction will be divided into three cases according to the following plan.
We consider a short exact sequence $0\to K\to F\to N\to0$, where $K$ is the $A$-torsion submodule of $F$ and $N$ is $A$-torsion-free.
Then $N$ can be treated by Auslander's techniques (as extended in Section~\ref{sec:almost-fin-modules}); see Case~(1) of the proof below. The $A$-torsion module $K$ is supported over a subgerm of $X_0$ of strictly smaller fibre dimension over $Y_0$, so it can be treated by induction (Case~(2)). For the general case (3), we want to show that the analytic tensor powers of either $N$ or $K$ \emph{embed} into the corresponding powers of $F$, and hence so do their $R$-vertical elements.

The latter would be automatic if $F$ were the \emph{direct} sum of $K$ and $N$. This is not true, in general.
We can, however, choose $N$ to be a submodule of $F$ of the form $g\!\cdot\!F$, for a suitable choice of $g\in A$.
It follows that the analytic tensor powers of $K$ embed into those of $F$, unless $gF$ is not 
$R$-flat. In the latter case, in turn, we show that the analytic tensor powers of $gF$ 
embed into those of $F$.
\par

The following lemma will be used in Case~(3) of the proof of Theorem~\ref{thm:main} below.
It will allow us to conclude that the analytic tensor powers of $gF$
embed into the corresponding powers of $F$, provided $gF$ is $A$-torsion-free.

\begin{lemma}
\label{lem:non-zerodiv}
Let $R=\C\{y_1,\dots,y_n\}$, and let $A$ and $B$ be regular local analytic $R$-algebras.
Suppose that $M$ and $N$ are finite $A$- and $B$-modules (respectively).
Let $g\in A$, $h\in B$, and $m\in gM\tensR hN$ all be nonzero elements.
If $m=0$ as an element of $M\tensR N$, then $(g\tensR h)\!\cdot\!m=0$ in $gM\tensR hN$.
In other words, if $g\tensR h$ is not a zero-divisor of $gM\tensR hN$, then the canonical
homomorphism $gM\tensR hN \to M\tensR N$ is an embedding.
\end{lemma}

\begin{proof}
Using the identification
\[
gM\tensR hN \cong (gM\tensR B)\otimes_{A\tensR B}(A\tensR hN)\,,
\]
we can write $m=\sum_{i=1}^km_i\otimes n_i$, where the $m_i\in gM\tensR B$, 
and $n_1,\dots,n_k$ generate $A\tensR hN$. The latter can be extended to a 
sequence $n_1,\dots,n_k,n_{k+1},\dots,n_t$ generating $A\tensR N$. 
Setting $m_{k+1}=\dots=m_t=0$, we get 
$m=\sum_{i=1}^tm_i\otimes n_i\in(M\tensR B)\otimes_{A\tensR B}(A\tensR N)$. 
By \cite[Lemma~6.4]{Eis}, $m=0$ in $M\tensR N$ if and only if there are 
$m'_1,\dots,m'_s\in M\tensR B$ and $a_{ij}\in A\tensR B$, such that
\begin{align}
\label{eqn12}
\sum_{j=1}^s a_{ij}m'_j &= m_i\ \ \mathrm{in}\,\ M\tensR B, 
\quad \mathrm{for\ all}\ i\,;\\
\label{eqn13}
\sum_{i=1}^t a_{ij}n_i &= 0\ \ \mathrm{in}\,\ A\tensR N, 
\quad \mathrm{for\ all}\ j\,.
\end{align}
Multiplying the equations (\ref{eqn12}) by $g\tensR1$, we get
\begin{equation}
\label{eqn14}
\sum_{j=1}^s a_{ij}(g\tensR1)m'_j=(g\tensR1)m_i\ \ \mathrm{in}\,\ gM\tensR B,
\quad \mathrm{for\ all}\ i\,;
\end{equation}
hence $(g\tensR1)m=0$ in $gM\tensR N$, by (\ref{eqn13}) and (\ref{eqn14}). 

Now write $(g\tensR1)m=\sum_{i=1}^lm_i\otimes n_i$, where the $n_i\in A\tensR N$, 
and $m_1,\dots,m_l$ generate $gM\tensR B$. Then $(g\tensR1)m=0$ in 
$(gM\tensR B)\otimes_{A\tensR B}(A\tensR N)$ if and only if there are 
$n'_1,\dots,n'_p\in A\tensR N$ and $b_{ij}\in A\tensR B$, such that
\begin{align}
\label{eqn15}
\sum_{j=1}^p b_{ij}n'_j &= n_i\ \ \mathrm{in}\,\ A\tensR N,
\quad \mathrm{for\ all}\ i\,;\\
\label{eqn16}
\sum_{i=1}^l b_{ij}m_i &= 0\ \ \mathrm{in}\,\ gM\tensR B,
\quad \mathrm{for\ all}\ j\,.
\end{align}
Multiplying the equations (\ref{eqn15}) by $1\tensR h$, we get
\begin{equation}
\label{eqn17}
\sum_{j=1}^p b_{ij}(1\tensR h)n'_j=(1\tensR h)n_i\ \ \mathrm{in}\,\ A\tensR hN,
\quad \mathrm{for\ all}\ i\,;
\end{equation}
hence $(g\tensR h)m=0$ in $gM\tensR hN$, by (\ref{eqn16}) and (\ref{eqn17}).
Thus $g\tensR h$ is a zero-divisor of $gM\tensR hN$, as required.
\end{proof}

\subsection{Proof of Theorem~\ref{thm:main}}
\label{subsec:proof}

Let $F$ be an almost finitely generated module over $R:=\C\{y_1,\dots,y_n\}$.
By Remark~\ref{rem:vert-element}, there exists $m\geq0$ such that $F$ is 
finitely generated as a module over $A:=R\{x\} = R\{x_1,\dots,x_m\}$.
Let $X$ and $Y$ be connected open neighbourhoods of the origins in $\C^{m+n}$ and $\C^n$ (respectively),
and let $\vp:X\to Y$ be the canonical coordinate projection.
Let $\FF$ be a coherent sheaf of $\OO_X$-modules whose stalk at the origin in $X$ equals $F$,
and let $\GG$ be a coherent $\OO_{\Xpn}$-module whose stalk at the origin $0^{\{n\}}$ in $\Xpn$ equals $\FpnR$.
We can identify $R$ with $\OO_{Y,0}$ and $A$ with $\OO_{X,0}$. Then $F$ is
$R$-flat if and only if $\FF_0$ is $\OO_{Y,0}$-flat.
\par

We first prove the ``only if'' direction of Theorem \ref{thm:main}, 
by contradiction. Assume that $F$ is $R$-flat. Since flatness is an open condition, by Douady's 
theorem \cite{Dou}, we can assume that $\FF$ and $\GG$ are 
$\OO_Y$-flat. Suppose that $\FpnR$ has a vertical element over $\OO_{Y,0}$.
In other words, (after shrinking $X$ and $Y$ if necessary) there exist a nonzero 
section $\tilde{m}\in\GG$ and an analytic subset $Z\subset\Xpn$, such that 
$Z_0=\VV(\Ann_{\OO_{\Xpn,0^{\{n\}}}}(\tilde{m}_0))$ and the image 
$\vpn(Z)$ has empty interior in $Y$. Let $\tilde{\vp}$ denote the restriction
$\vpn|_Z:Z\to Y$. Consider $\xi\in Z$ such that the fibre dimension of 
$\tilde{\vp}$ at $\xi$ is minimal. Then the fibre dimension $\fbd_x\tilde{\vp}$ 
is constant on some open neighbourhood $U$ of $\xi$ in $Z$. By the Remmert Rank 
Theorem, $\tilde{\vp}(U)$ is locally analytic in $Y$ near $\eta=\tilde{\vp}(\xi)$.
Since $\tilde{\vp}(Z)$ has empty interior in $Y$, it follows that there is a 
holomorphic function $g$ in a neighbourhood of $\eta$ in $Y$, such that 
$(\tilde{\vp}(U))_{\eta}\subset\VV(g_{\eta})$. Therefore,
$\tilde{\vp}^*_{\xi}(g_{\eta})\!\cdot\!\tilde{m}_{\xi}=0$ in $\GG_{\xi}$; 
i.e., $\GG_{\xi}$ has a (nonzero) zero-divisor in $\OO_{Y,\eta}$, contradicting flatness.
\par

We will now prove the more difficult ``if'' direction of the theorem, by induction on $m$.
If $m=0$, then $F$ is finitely generated over $R$, and the result follows from 
Auslander's theorem \ref{thm:Auslander} (because flatness of finitely generated
modules over a local ring is equivalent to freeness, the analytic tensor product 
equals the ordinary tensor product for finite modules, and vertical elements in 
finite modules are just zero-divisors).\par

The inductive step will be divided into three cases:
\begin{enumerate}
\item $F$ is torsion-free over $A$;
\item $F$ is a torsion $A$-module;
\item $F$ is neither $A$-torsion-free nor a torsion $A$-module.
\end{enumerate}
\par

\subsection*{Case (1)}
We prove this case independently of the inductive hypothesis. We essentially 
repeat the argument of Galligo and Kwieci{\'n}ski~\cite{GK}, which itself is 
an adaptation of Auslander \cite{Au} to the almost finitely generated context.
\par

Suppose that $\FpnR$ has no vertical elements over $R$. Then it follows from 
Lemma~\ref{lem:fd-sum}(1) that $F^{\antens^i_R}$ has no vertical elements,
for $i = 1,\ldots,n$. By Lemma~\ref{lem:fd-sum}(3),
\[
\mathrm{fd}_R(\FpnR) = \mathrm{fd}_R(F)+\mathrm{fd}_R(F^{\antens^{n-1}_R}) = \dots = n\cdot\mathrm{fd}_R(F)\,.
\]
On the other hand, since $\FpnR$ has no vertical elements over $R$, it has no 
zero-divisors over $R$, so that $\mathrm{depth}_R(\FpnR)\geq1$. It follows from 
Proposition~\ref{rem:GK-facts}(2) that $\mathrm{fd}_R(\FpnR)<n$. Hence $n\cdot\mathrm{fd}_R(F)<n$.
This is possible only if $\mathrm{fd}_R(F)=0$; i.e., $F$ is $R$-flat.

\subsection*{Case (2)}
Suppose that $F$ is not $R$-flat and a torsion $A$-module. We will show that then $\FpnR$ contains vertical
elements over $R$. Let $I=\Ann_A(F)$. Since  every element of $F$ is annihilated by 
some nonzero element of $A$, and $F$ is finitely generated over $A$, then $I$ 
is a nonzero ideal in $A$. Put $B=A/I$; then $F$ is finitely generated over $B$. 
Let $I(0)$ denote the \emph{evaluation} of $I$ \emph{at} $y=0$ (i.e., $I(0)$ is the ideal generated by $I$
in $A(0) := A \tensR R/\mm_R \cong \C\{x_1,\dots,x_m\}$).\par

First suppose that $I(0)\neq(0)$. Then there exists $g\in I$ such that 
$g(0):=g(0,x)\neq0$, and $F$ is a finite $A/(g)A$-module. It follows that 
(after an appropriate linear change in the $x$-coordinates) $g$ is regular in 
$x_m$ and hence, by the Weierstrass Preparation Theorem, that  $F$ is finite 
over $R\{x_1,\dots,x_{m-1}\}$. Therefore, $\FpnR$ has a vertical element over $R$, 
by the inductive hypothesis.\par

On the other hand, suppose that $I(0)=(0)$; i.e., $I\subset\mm_RA$.
Then $B\tensR R/\mm_R = (A/I)\tensR R/\mm_R$ equals $\C\{x_1,\dots,x_m\}$.
Let $Z$ be a closed analytic subspace of $X$ such that $\OO_{Z,0}\cong B$, and 
let $\tilde{\vp}:=\vp|_Z$. It follows that the fibre $\tilde{\vp}^{-1}(0)$
equals $\C^m$. Of course, $m$ is not the generic fibre dimension of $\tilde{\vp}$ on any irreducible component of $Z$,
because otherwise all its fibres would equal $\C^m$, so we would
have $B=A$ and $I=(0)$, contrary to the choice of $I$. Therefore, by 
Proposition \ref{prop:openness}, there is an isolated algebraic vertical 
component in the $n$-fold fibred power of $\tilde{\vp}_0$; i.e., $B^{\antens^n_R}$ 
has a zero-divisor in $R$. But $\FpnR$ is a faithful 
$B^{\antens^n_R}$-module, so itself it has a zero-divisor (hence a vertical 
element) over $R$.

\subsection*{Case (3)}
Suppose that $F$ is not $R$-flat, $F$ has zero-divisors in $A$, but 
$\Ann_A(F)=(0)$. Let
\[
K:=\{f\in F:af=0\ \mathrm{for\ some\ nonzero}\ a\in A\};
\]
i.e., $K$ is the $A$-torsion submodule of $F$. Since $K$ is a submodule of a finitely 
generated module over a Noetherian ring, $K$ is finitely generated; say 
$K=\sum_{i=1}^sA\!\cdot\!f_i$. Take $a_i\in A\setminus\{0\}$ such that $a_if_i=0$,
and put $g=a_1\cdots a_s$. Then the sequence of $A$-modules
\begin{equation}
\label{eqn5}
0\to K\to F\stackrel{\cdot g}{\rightarrow} gF\to0
\end{equation}
is exact, and $gF$ is a torsion-free $A$-module.\par

First suppose that $gF$ is $R$-flat. Then by applying $\tensR K$ and $F \tensR$ to (\ref{eqn5}), we get short exact sequences
\begin{align*}
&0 \to K\tensR K \to F\tensR K \to gF\tensR K \to 0,\\
&0 \to F\tensR K \to F\tensR F \to F\tensR\, gF \to 0\,.
\end{align*}
So we have injections
\[
K\tensR K\hookrightarrow F\tensR K\hookrightarrow F\tensR F\,,
\]
and by induction, an injection $K^{\antens^i_R}\hookrightarrow F^{\antens^i_R}$,
for all $i\geq1$. In particular, $K^{\antens^n_R}$ is a submodule of $\FpnR$. 
Since $gF$ is $R$-flat and $F$ is not $R$-flat, it follows that $K$ is not 
$R$-flat. Therefore, by Case (2), $K^{\antens^n_R}$ (and hence $\FpnR$) 
has a vertical element over $R$.
\par

Now suppose that $gF$ is not $R$-flat. Then $\gFpnR$ has a vertical element over $R$, by Case~(1).
We will show that $\gFpnR$ embeds into $\FpnR$, and hence so do its vertical elements.
By Lemma~\ref{lem:non-zerodiv}, in order for $\gFpnR$ to embed into $\FpnR$,
it suffices to prove that $g^{\tensRn}$ is not a zero-divisor of $\gFpnR$.
\par

To simplify the notation, let $B$ denote the ring $A^{\tensRn}$, and let $h:=g^{\tensRn}\in B$.
Since $\gFpnR$ is a finite $B$-module, we can write $\gFpnR=B^q/M$, where $q\geq1$ and 
$M$ is a $B$-submodule of $B^q$.
Given $b\in B$, let $(M:b)$ denote the $B$-submodule of $B^q$ consisting of those elements $m\in B^q$ for which $b\!\cdot\!m\in M$.
Since
\[
(M:h)\ \subset\ (M:h^2)\ \subset\dots\subset\ (M:h^l)\ \subset\dots
\]
is an increasing sequence of submodules of a Noetherian module $B^q$, it stabilizes;
i.e., there exists $k\geq1$ such that $(M:h^{k+1})=(M:h^k)$.
In other words, there exists $k\geq1$ such that $h$ is not a zero-divisor in $h^k\!\cdot\!B^q/M$; i.e.,
$g^{\tensRn}$ is not a zero-divisor in $(g^{k+1}F)^{\tensRn}$.
\par

Observe though that (because $gF$ is $A$-torsion-free),
multiplication by $g$ induces an isomorphism $gF\to g^2F$ of $A$-modules,
and in general, $gF\cong g^lF$, for $l\geq1$. We thus have isomorphisms $(gF)^{\tensRn}\cong(g^lF)^{\tensRn}$ 
of $B$-modules, for $l\geq1$. In particular, for every $l\geq1$, 
$g^{\tensRn}$ is a zero-divisor of $(gF)^{\tensRn}$ if and only if it is a zero-divisor of 
$(g^lF)^{\tensRn}$. Therefore, by Lemma~\ref{lem:non-zerodiv}, we have an embedding 
$(gF)^{\antens^n_R}\hookrightarrow\FpnR$. This completes the proof of Theorem~\ref{thm:main}.

\bibliographystyle{amsplain}

\end{document}